\documentclass[reqno,fleqn,a4paper,12pt]{amsart}

\usepackage[numbers,sort&compress]{natbib}
\usepackage{amssymb,amsmath,mathrsfs,amsthm}
\usepackage{graphics,epsfig}
\usepackage{graphicx}
\usepackage{enumitem}
\usepackage[T1]{fontenc}
\usepackage{pgf,tikz}
\usetikzlibrary{arrows}
\theoremstyle{plain}
\newtheorem{thm}{Theorem}[section]
\newtheorem{cor}[thm]{Corollary}
\newtheorem{lem}[thm]{Lemma}

\theoremstyle{definition}
\newtheorem{dfn}[thm]{Definition}

\theoremstyle{remark}
\newtheorem{rmk}[thm]{Remark}

\newcommand{\reDeclareMathOperator}[2]{\let#1\undefined \DeclareMathOperator{#1}{#2}}

\reDeclareMathOperator{\mod}{mod}

\numberwithin{subsection}{section} \numberwithin{equation}{section}

\begin{document}

\title[Random fixed point theorems for Hardy-Rogers self-random operators]{Random fixed point theorems for  Hardy-Rogers self-random operators with applications to random integral equations}

\author[P. Saipara]{Plern Saipara$^1$}

\author[P. Kumam]{Poom Kumam$^{1,2,3,*}$}

\author[Y. J. Cho]{Yeol Je Cho$^{3}$}

\address{$^1$KMUTTFixed Point Research Laboratory, Department of Mathematics, Room SCL 802 Fixed Point Laboratory, Science Laboratory Building, King Mongkut's University of Technology Thonburi (KMUTT), 126 Pracha-Uthit Road, Bang Mod, Thung Khru, Bangkok 10140 Thailand}
\address{$^{2}$KMUTT-Fixed Point Theory and  Applications Research Group, Theoretical and Computational Science Center (TaCS), Science Laboratory Building, Faculty of Science, King Mongkut's University of Technology Thonburi (KMUTT), 126 Pracha-Uthit Road, Bang Mod, Thung Khru, Bangkok 10140 Thailand.}

\address{$^{3}$Department of Medical Research, China Medical University Hospital, China Medical University, Taichung 40402, Taiwan}
\address{$^4$Department of Mathematics Education, Gyeongsang Natoinal University, Jinju 660-701, Korea, and Center for General Education, China Medical University,  Taichung, 40402, Taiwan}
\email{$^1$plern.spn@mail.kmutt.ac.th (Plern Saipara)}
\email{$^{1,2,3,*}$poom.kum@kmutt.ac.th (Poom Kumam) (Corresponding author)}
\email{$^4$yjcho@gnu.ac.kr (Yeol Je Cho)}

\thanks{$^*$Corresponding author: poom.kum@kmutt.ac.th (P. Kumam)\\
This project was supported by the Theoretical and Computational Science (TaCS) Center.}

\begin{abstract}
    In this paper, we prove some random fixed point theorems for  Hardy-Rogers self-random operators in separable Banach spaces and, as some applications, we show the existence of a solution for  random nonlinear integral equations in Banach spaces. Some stochastic versions of deterministic fixed point theorems for Hardy-Rogers self mappings and stochastic  integral equations  are obtained.
\end{abstract}
\maketitle

\section{{\bf Introduction}}
\vskip 3mm

 Some well known random fixed point theorems are stochastic generalizations of Banach's fixed point theorem and Banach's type fixed point theorems in complete metric spaces. In 1955, Spacek \cite{Spacek} and Hans  \cite{Hans1, Hans2} initiated to prove random fixed point theorems for random contraction mappings in separable complete metric spaces. In 1966, Mukherjee \cite{Mukherjee} proved a random fixed point theorem in the sense of Schaduer's fixed point theorem in atomic probability measure spaces. Especially, in 1976, the work of Bharucha-Reid \cite{BharuchaReid} has been developed by various mathematicians. In 1979, Itoh \cite{Itoh} extended some random fixed point theorems of Spacek and Hans to the setting of multi-valued contraction mappings and  applied random fixed point theorems to solve some random differential equations in Banach spaces. In 1984, Sehgal and Waters \cite{SehgalWaters} proved some random fixed point theorems including  classical results given by  Rothe \cite{Rothe}.

  Recently, Beg and Shahzad \cite{BegShahzad} showed the existence of  random common  fixed points and random coincidence points of a pair of compatible random multi-valued mappings in Polish spaces. Especially, Kumam et al. \cite{KumamKumam01, KumamPlubtieng04, KumamPlubtieng08, Kumam02, KumamKumam02} proved many random fixed point theorems for multi-valued nonexpansive nonself-mappings  satisfying the inwardness condition in Banach spaces (see \cite{KumamPlubtieng05}). Jung et al. \cite{JungChoKangLeeThakur} proved random fixed point theorems for a certain
class of mappings in banach spaces. Cho et al. \cite{ChoLiHuang} proved random Ishikawa iterative sequence with errors for approximating random fixed points. Likewise, Kumam and Plubtieng \cite{KumamPlubtieng01} showed the existence of a random coincidence point for a pair of reciprocally continuous and compatible single-valued and multi-valued mappings and Saha \cite{Saha}, Saha and Debnath \cite{SahaDebnath} established some random fixed point theorems in  separable Hilbert spaces and  separable Banach spaces, respectively. On the other hand, Padgett \cite{Padgett}, Achari \cite{Achari}, Saha and Dey \cite{SahaDey} applied  some random fixed point theorems  to show the existence of solutions of  random nonlinear integral equations in Banach spaces.

   Recently, Saha and Ganguly \cite{SahaGanguly} proved some random fixed point theorems for a class of contractive mappings in  separable Banach spaces equipped with a complete probability measure.
   \vskip 2mm

In fact, Banach's contraction principle (\cite{Banach}) is very important to show the existence of solutions of some nonlinear equations, differential and integral equations, and other nonlinear problems. Since Banach's contraction principle, many authors have studied in several ways.
\vskip 2mm

\begin{thm} {\rm (Banach's contraction principle)} If $(X, d)$ is a complete metric space and $T : X \rightarrow X$ be a mapping such that,
for some $\alpha \in [0, 1)$,
\begin{eqnarray} \label{IntroHardyRogerscontraction001}
d(Tx, Ty) \leq \alpha d(x, y)
\end{eqnarray}
for each $x, y \in X$, then $T$ has a unique fixed point in $X$.
\end{thm}
\vskip 2mm

Note that the mapping $T$ satisfying the Banach contraction condition is continuous, but the mappings $T$ satisfying the following contractions conditions are not continuous.
\vskip 2mm

(1)\, In 1968, Kannan's contraction (\cite{Rhoades}): for some $\beta \in [0, \frac{1}{2})$,
\begin{eqnarray} \label{IntroHardyRogerscontraction002}
d(Tx, Ty) \leq \beta [d(x, Tx) + d(y, Ty)]
\end{eqnarray}
for each $x, y \in X$;
\vskip 1mm

(2)\, In 1971, Reich's contraction (\cite{Rhoades}): for some $\alpha, \beta, \gamma \geq 0$ with $\alpha + \beta + \gamma < 1$,
\begin{eqnarray} \label{IntroHardyRogerscontraction003}
d(Tx, Ty) \leq \alpha d(x, Tx) + \beta d(y, Ty) + \gamma d(y, Ty)
\end{eqnarray}
for each $x, y \in X$;
\vskip 1mm

(3)\, In 1971, $\acute{\text{C}}$iri$\acute{\text{c}}$'s contraction (\cite{Rhoades}): for some $\alpha, \beta, \gamma, \delta \geq 0$ with $\alpha + \beta + \gamma + 2\delta < 1$,
\begin{eqnarray} \label{IntroHardyRogerscontraction004}
d(Tx, Ty) \leq \alpha d(x, y) + \beta d(x, Tx) + \gamma d(y, Ty) + \delta [d(x, Ty) + d(y, Tx)]
\end{eqnarray}
for each $x, y \in X$;
\vskip 1mm

(4)\, In 1972, Chatterjea's contraction (\cite{Rhoades}): for some $\beta \in [0, \frac{1}{2})$,
\begin{eqnarray} \label{IntroHardyRogerscontraction005}
d(Tx, Ty) \leq \beta [d(x, Ty) + d(y, Tx)]
\end{eqnarray}
for each $x, y \in X$;
\vskip 1mm

(5)\, In 1972, Zamfirescu contractive conditions (\cite{Rhoades}): there exist real numbers $\alpha, \beta, \gamma, 0 \leq \alpha < 1, 0 \leq \beta < 1, \gamma < \frac{1}{2}$, such that, for each $x, y \in X$, at least one of the following is true:
\begin{enumerate}
  \item [(i)] $d(Tx, Ty) \leq \alpha d(x, y)$;
  \item [(ii)] $d(Tx, Ty) \leq \beta [d(x, Tx) + d(y, Ty)]$;
  \item [(iii)] $d(Tx, Ty) \leq \gamma [d(x, Ty) + d(y, Tx)]$.
\end{enumerate}
For each $x, y \in X, x \neq y$,
\begin{eqnarray} \label{IntroHardyRogerscontraction005b}
d(Tx, Ty) < \max \{ d(x, y), [d(x, Tx) + d(y, Ty)]/2, [d(x, Ty) + d(y, Tx)]/2 \}.
\end{eqnarray}
\vskip 1mm

(6)\, In 1973, Hardy and Rogers's contraction (\cite{Rhoades}): for some $\alpha, \beta, \gamma, \delta, \eta \geq 0$ with $\alpha + \beta + \gamma + \delta + \eta < 1$,
\begin{eqnarray} \label{IntroHardyRogerscontraction006}
d(Tx, Ty) \leq \alpha d(x, y) + \beta d(x, Tx) + \gamma d(y, Ty) + \delta d(x, Ty) + \eta d(y, Tx)
\end{eqnarray}
for each $x, y \in X$:
\vskip 2mm

In 2000, $\acute{\text{C}}$iri$\acute{\text{c}}$ \cite{Ciric01} dealt with a class of mappings (not necessarily continuous) satisfying Gregus type contraction in metric spaces (\cite{Gregus01}) and proved the following fixed point theorem:
\vskip 2mm

\begin{thm}
Let $C$ be a closed convex subset of a complete convex metric space $X$ and $T : C \rightarrow C$ be a mapping satisfying
\begin{eqnarray} \label{IntroHardyRogerscontraction007}
  d(Tx, Ty) \leq ad(x, y) + b\max \{d(x, Tx), d(y, Ty)\}  + c [d(x, Ty) + d(y, Tx)]
\end{eqnarray}
for all $x, y \in C$, where $0 < a < 1$, $a + b = 1$ and $c \leq \frac{4 - a}{8 - b}$.  Then $T$ has a unique fixed point in $X$.
\end{thm}
\vskip 2mm

Moreover, Common fixed points under contractive conditions in cone metric spaces was studied by Radenovi$\acute{\text{c}}$ (see in \cite{Radenovic}).

Recently, Saha and Ganguly \cite{SahaGanguly} proved some random fixed point theorems  for a certain class of contractive mappings  in a separable Banach space equipped with a complete probability measure as follows:
\vskip 2mm

\begin{thm} \label{MainMotivate}
Let $X$ be a separable Banach space and $(\Omega, \beta, \mu)$ be a complete probability measure space. Let $T : \Omega \times X \rightarrow X$ be a continuous random operator such that
 for all $\omega \in \Omega$, $T$ satisfies
\begin{eqnarray} \label{IntroHardyRogerscontraction008}
 && \|T(\omega, x_1)-T(\omega, x_2)\|\nonumber\\
  &\leq& a(\omega) \max \{\|x_1 - x_2\|, \nonumber\\
  && \qquad\quad \qquad c(\omega) [\|x_1 - T(\omega, x_2)\| + \|x_2 - T(\omega, x_1)\|]\}\nonumber\\
    &&+ b(\omega) \max\{\|x_1 - T(\omega, x_1)\|,  \|x_2 - T(\omega, x_2)\|\} \\ \nonumber
\end{eqnarray}
for all random variables $x_1, x_2 \in X$ where $a(\omega), b(\omega), c(\omega)$ are real-valued random variables such that $0 < a(\omega) < 1$, $a(\omega) + b(\omega) = 1$, $c(\omega) \leq \frac{4 - a(\omega)}{8 - a(\omega)}$ almost surely. Then there exist unique random fixed point of $T$ in $X$.
\end{thm}
\vskip 3mm

Note that, if $\beta = 0$ or $\gamma = 0$ and $\delta = \eta$, then  fixed point theorems for Hardy and Roger's contraction \eqref{IntroHardyRogerscontraction006} reduced to  fixed point theorems for Gregus type contraction \eqref{IntroHardyRogerscontraction007}.
\vskip 2mm

The purpose of this paper is to prove some random fixed point theorems for  random Hardy-Rogers self-mappings in separable Banach spaces and, by using our main results, we show the existence of solutions of  random nonlinear integral equations.
\vskip 4mm

\section{{\bf Preliminaries}}
\vskip 3mm

Throughout this paper, X will denote a separable Banach over the real. Let $\beta_X$ be a $\sigma$-algebra of Borel subsets of $X$. Let $(\Omega, \beta, \mu)$ denote a complete probability measure space with the measure $\mu$ and $\beta$ be a
$\sigma$-algebra of subsets of $\Omega$. For more details,  see  Joshi and Bose \cite{JoshiBose}.
\vskip 2mm

\begin{dfn} \label{def2.1}
(1)\, A mapping $x : \Omega \rightarrow X$ is called  an {\it $X$-valued random variable} if the inverse image under the mapping $x$ of every Borel set $B$ of $X$ belongs to $\beta$, that is, $x^{-1}(X) \in \beta$ for all $B \in \beta_X$.

(2) A mapping $x : \Omega \rightarrow X$ is called a {\it finitely-valued random variable} if it is constant on each finite number of disjoint sets $A_i \in \beta$ and is equal to $0$ on $\Omega - (\bigcup_{i=1}^n A_i)$. The mapping $x$ is called a {\it simple random variable} if it is finitely valued and $\mu\{\omega : \|x(\omega)\| > 0\} < \infty$.

(3)\, A mapping $x : \Omega  \rightarrow X$ is called  a {\it strong random variable} if there exists a sequence $\{x_n(\omega)\}$ of simple random variables which converges to $x(\omega)$ almost surely, that is, there exists a set $A_0 \in \beta$ with $\mu(A_0) = 0$ such that $\lim_{n \rightarrow \infty} x_n(\omega) = x(\omega)$ for any $\omega \in \Omega - A_0$.

(4)\, A mapping $x : \Omega \rightarrow X$ is called  a {\it weak random variable} if the function $x^*(x(\omega))$ is a real-valued random variable for each $x^* \in X^*$, where $X^*$ denots the first normed dual space of $X$.
\end{dfn}
\vskip 2mm

In a separable Banach space $X$, the notions of strong and weak random variables $x : \Omega \rightarrow X$ (\cite{JoshiBose}) coincide and, in $X$, $x$ is termed as a random variable.
\vskip 2mm

Now, we recall the following:
\vskip 2mm

\begin{thm} \label{thm2.5}
{\rm (\cite{JoshiBose})} Let $x, y : \Omega \rightarrow X$ be strong random variables and $\alpha, \beta$ be constants. Then the following statements hold:
\vskip 1mm

{\rm (1)}\, $\alpha x(\omega) + \beta y(\omega)$ is a strong random variable.

{\rm (2)}\, If $f(\omega)$ is a real-valued random variable and $x(\omega)$ is a strong random variable, then
 $f(\omega) x(\omega)$ is a strong random variable.

{\rm (3)}\, If $x_n(\omega)$ is a sequence of strong random variables converging strongly to $x(\omega)$ almost
surely, that is, if there exists a set $A_0 \in \beta$ with $\mu(A_0) = 0$ such that
$$
\lim_{n \rightarrow \infty} \|x_n(\omega) - x(\omega)\| = 0
$$
 for any $\omega \notin A_0$, then $x(\omega)$ is a strong random variable.
\end{thm}
\vskip 2mm

\begin{rmk}
If $X$ is a separable Banach space, then every strong and also weak random variable is measurable in the sense of Definition \ref{def2.1}.
\end{rmk}
\vskip 2mm

Let $Y$ be an another Banach space. We also need the following definitions (see Joshi
and Bose \cite{JoshiBose}).
\vskip 2mm

\begin{dfn}
(1)\, A mapping $F : \Omega \times X \rightarrow Y$ is called a {\it random mapping} if $F(\omega, x) = Y(\omega)$ is a $Y$-valued random variable for all $x \in X$.

(2)\, A mapping $F : \Omega \times X \rightarrow Y$ is called  a {\it continuous random mapping} if the set of all $\omega \in \Omega$ for which $F(\omega, x)$ is a continuous function of $x$ has measure one.

(3)\, A mapping $F : \Omega \times X \rightarrow Y$ is said to be {\it demicontinuous} at the $x \in X$ if $\|x_n - x\| \rightarrow 0$ implies $F(\omega, x_n) \rightharpoonup F(\omega, x)$ almost surely.
\end{dfn}
\vskip 2mm

\begin{thm}
{\rm (\cite{JoshiBose})} Let $F : \Omega \times X \rightarrow Y$ be a demicontinuous random mapping where a Banach space $Y$ is separable. Then, for any $X$-valued random variable $x$, the function $F(\omega, x(\omega))$ is a $Y$-valued random variable.
\end{thm}
\vskip 2mm

\begin{rmk}
(\cite{JoshiBose}) Since a continuous random mapping is a demicontinuous random mapping, Theorem \ref{thm2.5} is also true for a continuous random mapping.
\end{rmk}
\vskip 2mm

Also, we recall the following definitions (see Joshi and Bose \cite{JoshiBose}):
\vskip 2mm

\begin{dfn}
(1)\, An equation of the type $F(\omega, x(\omega)) = x(\omega)$, where $F : \Omega \times X \rightarrow X$ is a random mapping, is called a {\it random fixed point equation}.

(2)\, Any mapping $x : \Omega \rightarrow X$ which satisfies the random fixed point equation $F(\omega, x(\omega)) = x(\omega)$ almost surely is called
 a {\it wide sense solution} of the fixed point equation.

 (3)\, Any $X$-valued random variable $x(\omega)$ which satisfies
 \begin{center}
 $\mu\{\omega : F(\omega, x(\omega)) = x(\omega)\} = 1$
 \end{center}
 is called a {\it random solution} of the fixed point equation or a random fixed point of $F$.
\end{dfn}
\vskip 2mm

\begin{rmk}
A random solution is a wide sense solution of the fixed point equation. But the converse is not necessarily true. This is evident from an example, under Remark 1, in Joshi and Bose \cite{JoshiBose}.
\end{rmk}
\vskip 4mm

\section{{\bf The main results}}
\vskip 3mm

Motivated and inspired by Theorem \ref{MainMotivate}, we proposed the definition as follows:
\vskip 2mm

\begin{dfn}
Let $T : \Omega \times X \rightarrow X$ be a continuous random mapping. The random mapping $T$ is called {\it Hardy-Rogers' contraction} if, for any $\omega \in \Omega$,
\begin{eqnarray} \label{HardyRogerscontraction001}
 && \|T(\omega, x_1(\omega))-T(\omega, x_2(\omega))\|\nonumber\\
  &\leq& \alpha_1(\omega) \|x_1(\omega) - x_2(\omega)\| + \alpha_2(\omega) \|x_1(\omega) - T(\omega, x_1(\omega))\| \\ \nonumber
                        && + \alpha_3(\omega) \|x_2(\omega) - T(\omega, x_2(\omega))\| + \alpha_4(\omega) \|x_1(\omega) - T(\omega, x_2(\omega))\| \\ \nonumber
                        && + \alpha_5(\omega) \|x_2(\omega) - T(\omega, x_1(\omega))\|
\end{eqnarray}
for all random variables $x_1, x_2 : \Omega \rightarrow X$ and $\alpha_i : \Omega \rightarrow \mathbb{R}_+ \cup \{0\}$ for $i = 1, 2, 3, 4, 5$ such that $\sum^5_{i=1} \alpha_i(\omega) < 1$.
\end{dfn}
\vskip 2mm

\begin{thm} \label{Maintheorem1}
Let $X$ be a separable Banach space and $(\Omega, \beta, \mu)$ be a complete probability measure space. Let $T : \Omega \times X \rightarrow X$ be a continuous random mapping
satisfying Hardy-Rogers' contraction. Then there exists a unique random fixed point of $T$ in $X$.
\end{thm}
\vskip 2mm

\begin{proof}
Let
\begin{eqnarray*}
A = \{\omega \in \Omega : T(\omega, x_1) ~~~~ \text{is a continuous function of} ~~~~ x \},
\end{eqnarray*}
\begin{eqnarray*}
B = \Big\{\omega \in \Omega : \sum^5_{i=1} \alpha_i(\omega) < 1 \Big\}
\end{eqnarray*}
and
\begin{eqnarray*}
C_{x_1, x_2} &=& \{\omega \in \Omega : \|T(\omega, x_1(\omega))-T(\omega, x_2(\omega))\| \leq \alpha_1(\omega) \|x_1(\omega) - x_2(\omega)\| \\
                &&\qquad\quad + \alpha_2(\omega) \|x_1(\omega) - T(\omega, x_1(\omega))\| + \alpha_3(\omega) \|x_2(\omega) - T(\omega, x_2(\omega))\| \\
                &&\qquad\quad + \alpha_4(\omega) \|x_1(\omega) - T(\omega, x_2(\omega))\| + \alpha_5(\omega) \|x_2(\omega) - T(\omega, x_1(\omega))\| \}
\end{eqnarray*}
Let $S$ be a countable dense subset of $X$. Now, we prove that
\begin{center}
  $\bigcap_{x_1, x_2 \in X}(C_{x_1, x_2} \cap A \cap B) = \bigcap_{s_1, s_2 \in S}(C_{s_1, s_2} \cap A \cap B)$.
\end{center}
Now, for all $s_1, s_2 \in S$, we have
\begin{eqnarray} \label{HardyRogerscontraction002}
&&  \|T(\omega, s_1(\omega))-T(\omega, s_2(\omega))\|\nonumber\\
 &\leq& \alpha_1(\omega) \|s_1(\omega) - s_2(\omega)\| + \alpha_2(\omega) \|s_1(\omega) - T(\omega, s_1(\omega))\| \\ \nonumber
                    && + \alpha_3(\omega) \|s_2(\omega) - T(\omega, s_2(\omega))\| + \alpha_4(\omega) \|s_1(\omega) - T(\omega, s_2(\omega))\| \\ \nonumber
                    && + \alpha_5(\omega) \|s_2(\omega) - T(\omega, s_1(\omega))\|.
\end{eqnarray}
Since $S$ is dense in $X$, for any $\delta_i(x_i) > 0$, there exist $s_1, s_2 \in S$ such that $\|x_i - s_i\| < \delta_i(x_i)$ for each $i = 1, 2$.
Note that, for any  $x_1, x_2 \in X$,
\begin{eqnarray} \label{HardyRogerscontraction003}
  \|s_1(\omega) - s_2(\omega)\| \leq \|s_1(\omega) - x_1(\omega)\| + \|x_1(\omega) - x_2(\omega)\| + \|x_2(\omega) - s_2(\omega)\|,
\end{eqnarray}
\begin{eqnarray} \label{HardyRogerscontraction004}
  \|s_1(\omega) - T(\omega, s_1(\omega))\| &\leq& \|s_1(\omega) - x_1(\omega)\| + \|x_1(\omega) - T(\omega, x_1(\omega)\| \\ \nonumber
            && + \|T(\omega, x_1(\omega) - T(\omega, s_1(\omega)\|,
\end{eqnarray}
\begin{eqnarray} \label{HardyRogerscontraction005}
  \|s_2(\omega) - T(\omega, s_2(\omega))\| &\leq& \|s_2(\omega) - x_2(\omega)\| + \|x_2(\omega) - T(\omega, x_2(\omega)\| \\ \nonumber
            && + \|T(\omega, x_2(\omega) - T(\omega, s_2(\omega)\|,
\end{eqnarray}
\begin{eqnarray} \label{HardyRogerscontraction006}
  \|s_1(\omega) - T(\omega, s_2(\omega))\| &\leq& \|s_1(\omega) - x_1(\omega)\| + \|x_1(\omega) - T(\omega, x_2(\omega)\| \\ \nonumber
            && + \|T(\omega, x_2(\omega) - T(\omega, s_2(\omega)\|
\end{eqnarray}
and
\begin{eqnarray} \label{HardyRogerscontraction007}
  \|s_2(\omega) - T(\omega, s_1(\omega))\| &\leq& \|s_2(\omega) - x_2(\omega)\| + \|x_2(\omega) - T(\omega, x_1(\omega)\| \\ \nonumber
            && + \|T(\omega, x_1(\omega) - T(\omega, s_1(\omega)\|.
\end{eqnarray}
Suppose that
\begin{eqnarray} \label{HardyRogerscontraction008}
 && \|T(\omega, s_1(\omega))-T(\omega, s_2(\omega))\|\nonumber\\
  &\leq& \alpha_1(\omega) \|s_1(\omega) - s_2(\omega)\| + \alpha_2(\omega) \|s_1(\omega) - T(\omega, s_1(\omega))\| \\ \nonumber
            && + \alpha_3(\omega) \|s_2(\omega) - T(\omega, s_2(\omega))\|  + \alpha_4(\omega) \|s_1(\omega) - T(\omega, s_2(\omega))\| \\ \nonumber
            && + \alpha_5(\omega) \|s_2(\omega) - T(\omega, s_1(\omega))\|.
\end{eqnarray}
Since
\begin{eqnarray} \label{HardyRogerscontraction009}
  &&\|T(\omega, x_1(\omega))-T(\omega, x_2(\omega))\|\nonumber\\
   &\leq& \|T(\omega, x_1(\omega))-T(\omega, s_1(\omega))\|  + \|T(\omega, s_1(\omega))-T(\omega, s_2(\omega))\| \\ \nonumber
            && + \|T(\omega, s_2(\omega))-T(\omega, x_2(\omega))\|,
\end{eqnarray}
substituting \eqref{HardyRogerscontraction008} in \eqref{HardyRogerscontraction009}, we have
\begin{eqnarray} \label{HardyRogerscontraction010}
&&\|T(\omega, x_1(\omega))-T(\omega, x_2(\omega))\|\nonumber\\
 &\leq& \|T(\omega, x_1(\omega))-T(\omega, s_1(\omega))\| + \|T(\omega, s_2(\omega))-T(\omega, x_2(\omega))\| \\ \nonumber
            && + \alpha_1(\omega) \|s_1(\omega) - s_2(\omega)\|  + \alpha_2(\omega) \|s_1(\omega) - T(\omega, s_1(\omega))\| \\ \nonumber
            && + \alpha_3(\omega) \|s_2(\omega) - T(\omega, s_2(\omega))\|  + \alpha_4(\omega) \|s_1(\omega) - T(\omega, s_2(\omega))\| \\ \nonumber
            && + \alpha_5(\omega) \|s_2(\omega) - T(\omega, s_1(\omega))\|.
\end{eqnarray}
Thus, from \eqref{HardyRogerscontraction003}, \eqref{HardyRogerscontraction004}, \eqref{HardyRogerscontraction005}, \eqref{HardyRogerscontraction006}, \eqref{HardyRogerscontraction007}, \eqref{HardyRogerscontraction010}, it follows that
\begin{eqnarray} \label{HardyRogerscontraction011}
&&\|T(\omega, x_1(\omega))-T(\omega, x_2(\omega))\|\nonumber\\
 &\leq& \alpha_1(\omega)\|x_1(\omega)-x_2(\omega)\| + \alpha_2(\omega)\|x_1(\omega)- T(\omega, x_1(\omega))\| \\ \nonumber
            && + \alpha_3(\omega)\|x_2(\omega)- T(\omega, x_2(\omega))\|  + \alpha_4(\omega)\|x_1(\omega)- T(\omega, x_2(\omega))\| \\ \nonumber
            && + \alpha_5(\omega)\|x_2(\omega)- T(\omega, x_1(\omega))\| \\ \nonumber
            && + (1 + \alpha_2(\omega) + \alpha_5(\omega))\|T(\omega, x_1(\omega))-T(\omega, s_1(\omega))\| \\ \nonumber
            && +(1 + \alpha_3(\omega) + \alpha_4(\omega))\|T(\omega, x_2(\omega))-T(\omega, s_2(\omega))\| \\ \nonumber
            && +(\alpha_1(\omega) + \alpha_2(\omega) + \alpha_4(\omega))\|s_1(\omega)-x_1(\omega)\| \\ \nonumber
            && +(\alpha_1(\omega) + \alpha_2(\omega) + \alpha_4(\omega))\|s_2(\omega)-x_2(\omega)\| .
\end{eqnarray}
For any $\omega \in \Omega$, sice $T(\omega, x(\omega))$ is a continuous function of $x(\omega)$, for any $\varepsilon > 0$, there exists $\delta_i(x_i(\omega)) > 0 ~~~ (i = 1, 2)$ such that
\begin{eqnarray} \label{HardyRogerscontraction012}
\| T(\omega, x_1(\omega)) - T(\omega, s_1(\omega))\| < \frac{\varepsilon}{8}
\end{eqnarray}
whenever $\| x_1(\omega) - s_1(\omega) \| < \delta_1(x_1(\omega))$ and
\begin{eqnarray} \label{HardyRogerscontraction013}
\| T(\omega, x_2(\omega)) - T(\omega, s_2(\omega))\| < \frac{\varepsilon}{8}
\end{eqnarray}
whenever $\| x_2(\omega) - s_2(\omega) \| < \delta_1(x_2(\omega))$. Now, choosing
\begin{eqnarray} \label{HardyRogerscontraction014}
\delta_1 = \min\Big\{\delta_1(x_1(\omega)), \frac{\varepsilon}{8}\Big\}
\end{eqnarray}
and
\begin{eqnarray} \label{HardyRogerscontraction015}
\delta_2 = \min\Big\{\delta_2(x_2(\omega)), \frac{\varepsilon}{8}\Big\},
\end{eqnarray}
by \eqref{HardyRogerscontraction011}, we have
\begin{eqnarray*}
&&\|T(\omega, x_1(\omega))-T(\omega, x_2(\omega))\|\nonumber\\
 &\leq& \alpha_1(\omega)\|x_1(\omega)-x_2(\omega)\|  + \alpha_2(\omega)\|x_1(\omega)- T(\omega, x_1(\omega))\| \\ \nonumber
            && + \alpha_3(\omega)\|x_2(\omega)- T(\omega, x_2(\omega))\|  + \alpha_4(\omega)\|x_1(\omega)- T(\omega, x_2(\omega))\| \\ \nonumber
            && + \alpha_5(\omega)\|x_2(\omega)- T(\omega, x_1(\omega))\| + (1 + \alpha_2(\omega) + \alpha_5(\omega)) \frac{\varepsilon}{8} \\ \nonumber
            && + (1 + \alpha_3(\omega) + \alpha_4(\omega))\frac{\varepsilon}{8}  + (\alpha_1(\omega) + \alpha_2(\omega) + \alpha_4(\omega)) \frac{\varepsilon}{8} \\ \nonumber
            && + (\alpha_1(\omega) + \alpha_2(\omega) + \alpha_4(\omega))\frac{\varepsilon}{8}
\end{eqnarray*}
and so
\begin{eqnarray*}
&&\|T(\omega, x_1(\omega))-T(\omega, x_2(\omega))\|\nonumber\\
 &\leq& (2 + 2\sum^5_{i=1}\alpha_i(\omega))\frac{\varepsilon}{8}  + \alpha_1(\omega)\|x_1(\omega)-x_2(\omega)\| \\ \nonumber
            && + \alpha_2(\omega)\|x_1(\omega)- T(\omega, x_1(\omega))\| + \alpha_3(\omega)\|x_2(\omega)- T(\omega, x_2(\omega))\| \\ \nonumber
            && + \alpha_4(\omega)\|x_1(\omega)- T(\omega, x_2(\omega))\| + \alpha_5(\omega)\|x_2(\omega)- T(\omega, x_1(\omega))\|. \\ \nonumber
\end{eqnarray*}
Since $\varepsilon > 0$ is arbitrary, it follows that
\begin{eqnarray} \label{HardyRogerscontraction015}
&&\|T(\omega, x_1(\omega))-T(\omega, x_2(\omega))\|\nonumber\\
 &\leq& \alpha_1(\omega)\|x_1(\omega)-x_2(\omega)\|  + \alpha_2(\omega)\|x_1(\omega)- T(\omega, x_1(\omega))\| \\ \nonumber
            && + \alpha_3(\omega)\|x_2(\omega)- T(\omega, x_2(\omega))\|  + \alpha_4(\omega)\|x_1(\omega)- T(\omega, x_2(\omega))\| \\ \nonumber
            && + \alpha_5(\omega)\|x_2(\omega)- T(\omega, x_1(\omega))\|.
\end{eqnarray}
Thus we have $\omega \in \bigcap_{x_1, x_2 \in X}(C_{x_1, x_2} \cap A \cap B)$, which implies that
\begin{center}
$\bigcap_{s_1, s_2 \in S}(C_{s_1, s_2} \cap A \cap B) \subset \bigcap_{x_1, x_2 \in X}(C_{x_1, x_2} \cap A \cap B)$.
\end{center}
Also, we have
\begin{center}
$\bigcap_{x_1, x_2 \in X}(C_{x_1, x_2} \cap A \cap B) \subset \bigcap_{s_1, s_2 \in S}(C_{s_1, s_2} \cap A  \cap B)$.
\end{center}
Therefore, we have
\begin{center}
$\bigcap_{s_1, s_2 \in S}(C_{s_1, s_2} \cap A  \cap B) = \bigcap_{x_1, x_2 \in X}(C_{x_1, x_2} \cap A  \cap B)$.
\end{center}
Let $N' = \bigcap_{s_1, s_2 \in S}(C_{s_1, s_2} \cap A  \cap B)$. Then $\mu(N') = 1$, which implies that $T(\omega, x)$ is a deterministic mapping.
Hence $T$ has a unique random fixed point in $X$. This completes the proof.
\end{proof}
\vskip 2mm

If $\alpha_4(\omega) = \alpha_5(\omega) = 0$ in Theorem \ref{Maintheorem1}, then we obtain the following random fixed point theorem for Reich's contraction:
\vskip 2mm

\begin{cor} \label{RandomRiech}
Let $X$ be a separable Banach space and $(\Omega, \beta, \mu)$ be a complete probability measure space. Let $T : \Omega \times X \rightarrow X$ be a continuous random mapping
satisfying the following condition: for any $\omega \in \Omega$,
\begin{eqnarray*}
 && \|T(\omega, x_1(\omega))-T(\omega, x_2(\omega))\| \nonumber\\
 &\leq& \alpha_1(\omega) \|x_1(\omega) - x_2(\omega)\|  + \alpha_2(\omega) \|x_1(\omega) - T(\omega, x_1(\omega))\| \\ \nonumber
                                          && + \alpha_3(\omega) \|x_2(\omega) - T(\omega, x_2(\omega))\| \\ \nonumber
\end{eqnarray*}
for all random variables $x_1, x_2 : \Omega \rightarrow X$ and $\alpha_i : \Omega \rightarrow \mathbb{R}_+ \cup \{0\}$ for $i = 1, 2, 3$ such that $\sum^3_{i=1} \alpha_i(\omega) < 1$. Then there exists a unique random fixed point of $T$ in $X$.
\end{cor}
\vskip 2mm

If $\alpha_1(\omega) = \alpha_4(\omega) = \alpha_5(\omega) = 0$ in Theorem \ref{Maintheorem1}, then we obtain the following random fixed point theorem for Kannan's contraction:
\vskip 2mm

\begin{cor} \label{Kannan}
Let $X$ be a separable Banach space and $(\Omega, \beta, \mu)$ be a complete probability measure space. Let $T : \Omega \times X \rightarrow X$ be a continuous random mapping
satisfying the following condition:  for any $\omega \in \Omega$,
\begin{eqnarray*}
  &&\|T(\omega, x_1(\omega))-T(\omega, x_2(\omega))\|\nonumber\\
   &\leq& \alpha_2(\omega) \|x_1(\omega) - T(\omega, x_1(\omega))\| + \alpha_3(\omega) \|x_2(\omega) - T(\omega, x_2(\omega))\| \\ \nonumber
\end{eqnarray*}
for all random variables $x_1, x_2 : \Omega \rightarrow X$ and $\alpha_i : \Omega \rightarrow \mathbb{R}_+ \cup \{0\}$ for $i = 2, 3$ such that $\alpha_2(\omega) + \alpha_3(\omega) < 1$. Then there exists a unique random fixed point of $T$ in $X$.
\end{cor}
\vskip 2mm

If $\alpha_1(\omega) = \alpha_2(\omega) = \alpha_3(\omega) = 0$ in Theorem \ref{Maintheorem1}, then we obtain the following random fixed point theorem for Chatterjea's contraction:
\vskip 2mm

\begin{cor} \label{Chatterjea}
Let $X$ be a separable Banach space and $(\Omega, \beta, \mu)$ be a complete probability measure space. Let $T : \Omega \times X \rightarrow X$ be a continuous random mapping
satisfying the following condition:  for all $\omega \in \Omega$,
\begin{eqnarray*}
  &&\|T(\omega, x_1(\omega))-T(\omega, x_2(\omega))\|\nonumber\\
   &\leq& \alpha_4(\omega) \|x_1(\omega) - T(\omega, x_2(\omega))\| + \alpha_5(\omega) \|x_2(\omega) - T(\omega, x_1(\omega))\| \\ \nonumber
\end{eqnarray*}
for all random variables $x_1, x_2 : \Omega \rightarrow X$ and $\alpha_i : \Omega \rightarrow \mathbb{R}_+ \cup \{0\}$ for $i = 4, 5$ such that $\alpha_4(\omega) + \alpha_5(\omega) < 1$. Then there exists a unique random fixed point of $T$ in $X$.
\end{cor}
\vskip 2mm

\begin{rmk}
The random fixed point theorems for Hardy-Rogers's contraction reduced to the random fixed point theorems for $\acute{\text{C}}$iri$\acute{\text{c}}$'s contraction.
\end{rmk}
\vskip 4mm

\section{{\bf Applications to  random nonlinear integral equations}}
\vskip 3mm

In this section, we give an application of Theorem \ref{Maintheorem1} to show the existence and uniqueness of a solution of a nonlinear stochastic integral equation of the Hammerstein type (\cite{Padgett}):
$$
\aligned
x(t;\omega)=h(t;\omega)+\int_Sk(t;s;\omega)f(s;x(s;\omega))d\mu(s),
\endaligned
\eqno{(4.1)}
$$
where

(a)\, $S$ is a locally compact metric space with metric d defined on $S\times S$
and $\mu_0$ is a complete $\sigma$-finite measure defined on the collection of Borel
subsets of $S;$

(b)\, $\omega\in \Omega$ where $\omega$ is the supporting set of the probability measure
space $(\Omega, \beta, \mu);$

(c)\,  $x(t ; \omega)$ is the unknown vector-valued random variable for each $t \in S;$

(d)\, $h(t; \omega)$ is the stochastic free term defined for $t \in S;$

(e)\, $k(t, s; \omega)$ is the stochastic kernel defined for $t$ and $s$ in $S;$

(f)\, $f(t, x)$ is a vector-valued function of $t \in S$ and $x.$
\vskip 2mm

Note that the integral in the equation (4.1) is interpreted as a Bochner integral (\cite{Yosida}).
\vskip 2mm

Further, we assume that the union of a countable family $\{C_n\}$ of compact sets with $C_{n+1}\subset C_n$ is defined as $S$ such that,  for each other compact set in $S$,
  there exists $C_i$ which contains it (see \cite{Arens}).
\vskip 2mm

We define $C=C(S, L_2(\Omega, \beta, \mu))$ as a space of all continuous functions from $S$ into the space $L_2(\Omega, \beta, \mu)$ with the topology of uniform convergence on compact sets of $S$, that is, $x(t; \omega)$ is a vector-valued random variable for each fixed $t \in S$ such that
\begin{center}
$\| x(t; \omega) \|^2_{L_2(\Omega, \beta, \mu)} = \int_\Omega |x(t; \omega)|^2 d\mu(\omega) < \infty.$
\end{center}
Noted that $C(S, L_2(\Omega, \beta, \mu))$ is a space of locally convex (\cite{Yosida}) whose topology is defined by the countable family of semi-norms given by
\begin{center}
$\| x(t; \omega) \|_n = \sup_{t \in C_n} \| x(t; \omega) \|_{L_2(\Omega, \beta, \mu)}$
\end{center}
for each $n\geq 1$. Furthermore, since $L_2(\Omega, \beta, \mu)$ is complete,  $C(S, L_2(\Omega, \beta, \mu))$ is complete relative to this topology.
\vskip 2mm

 Next, we define $BC=BC(S, L_2(\Omega, \beta, \mu))$ as a Banach space of all bounded continuous functions from $S$ into $L_2(\Omega, \beta, \mu)$ with the norm
\begin{center}
$\| x(t; \omega) \|_{BC} = \sup_{t \in S} \| x(t; \omega) \|_{L_2(\Omega, \beta, \mu)}.$
\end{center}
The space $BC \subset C$ is a space of all second order vector-valued stochastic processes defined on $S$ which are bounded and continuous in mean-square.
 \vskip 2mm

Now, we consider the functions $h(t; \omega)$ and $f(t, x(t; \omega))$ to be in the $C(S, L_2(\Omega, \beta, \mu))$ space with respect to the stochastic kernel and  assume that,
 for each pair $(t, s)$, $k(t, s; \omega) \in L_\infty(\Omega, \beta, \mu)$ and the norm denoted by
\begin{center}
$\|| k(t, s; \omega)|\| = \| k(t, s; \omega)\|_{L_\infty(\Omega, \beta, \mu)} = \mu-ess \sup_{\omega \in \Omega}| k(t, s; \omega) |.$
\end{center}

Also, we suppose that $k(t, s; \omega)\in L_\infty(\Omega, \beta, \mu)$ is such that
$$
\|| k(t, s; \omega)|\| = \| x(s; \omega)\|_{L_2(\Omega, \beta, \mu)}
$$
is $\mu$-integrable with respect to $s$ for each $t \in S$ and $x(s; \omega) \in C(S, L_2(\Omega, \beta, \mu))$ and there exists a real-valued function $G$ $\mu$-a.e. on $S$
such that $G(S) \| x(s; \omega) \|_{L_2(\Omega, \beta, \mu))}$ is $\mu$-integrable and, for each pair $(t, s) \in S \times S$,
\begin{equation*}
\|| k(t, u; \omega) -  k(s, u; \omega)|\| \cdot \| x(u; \omega)\|_{L_2(\Omega, \beta, \mu)} \leq G(u)\| x(u; \omega)\|_{L_2(\Omega, \beta, \mu)} ~~~~ \mu-a.e.
\end{equation*}
Forward, assume that, for almost all $s \in S$, $k(t, s; \omega)$ is continuous in $t$ from $S$ into $L_\infty(\Omega, \beta, \mu)$.
\vskip 2mm

Now, we define the random integral operator $T$ on $C(S, L_2(\Omega, \beta, \mu))$ by
$$
\aligned
(Tx)(t; \omega) = \int_Sk(t, s; \omega)x(s; \omega)d\mu(s),
\endaligned
\eqno{(4.2)}
$$
where the integral is a Bochner integral. From the conditions on $k(t, s; \omega)$, it follows that, for each $t \in S$, $(Tx)(t; \omega) \in L_2(\Omega, \beta, \mu)$ and $(Tx)(t; \omega)$ is continuous in mean square by Lebesgue's dominated convergence theorem, that is, $(Tx)(t; \omega) \in C(S, L_2(\Omega, \beta, \mu))$.
\vskip 2mm

\begin{lem} {\rm (\cite{Padgett})}
The linear operator $T$ defined by the equation $(4.2)$ is continuous
from $C(S, L_2(\Omega, \beta, \mu))$ into itself.
\end{lem}
\vskip 2mm

\begin{proof}
See \cite{Padgett}.
\end{proof}
\vskip 2mm

\begin{dfn} (\cite{Achari}, \cite{LeePadgett})
Let $B$ and $D$ be Banach spaces. The pair $(B, D)$ is said to be \emph{admissible} with respect to a linear operator $T$ if $T(B) \subset D$.
\end{dfn}
\vskip 2mm

\begin{lem} {\rm (\cite{Padgett})}
If $T$ is a continuous linear operator from $C(S, L_2(\Omega, \beta, \mu))$ into itself and $B, D \subset C(S, L_2(\Omega, \beta, \mu))$ are Banach spaces stronger than $C(S, L_2(\Omega, \beta, \mu))$ such that $(B, D)$ is admissible with respect to $T$, then $T$ is continuous from $B$ into $D$.
\end{lem}
\vskip 2mm

By a {\it random solution} of the equation (4.1), we mean a function
$$
x(t; \omega) \in C(S, L_2(\Omega, \beta, \mu))
$$
which satisfies the equation (4.1) $\mu-a.e.$
\vskip 2mm

Now, by using Theorem \ref{Maintheorem1}, we prove the following:
\vskip 2mm

\begin{thm} \label{Maintheorem2}
If  the stochastic integral equation $(4.1)$ is subject to the following conditions:

 {\rm (1)}\, $B$ and $D$ are Banach spaces stronger than $C(S, L_2(\Omega, \beta, \mu))$ such that $(B,D)$ is admissible with respect to the integral operator defined by $(4.2)$;

   {\rm (2)}\, $x(t; \omega) \mapsto f(t, x(t; \omega))$ is an operator from the set $Q(\rho) = \{ x(t; \omega) : x(t; \omega) \in D, \| x(t; \omega) \|_D \leq \rho \}$ into the space $B$ satisfying
$$
\aligned
  &\quad\|f(t, x_1(t, \omega))-f(t, x_2(t, \omega))\|_B\\
   &\leq \alpha_1(\omega) \|x_1(t, \omega) - x_2(t, \omega)\|  + \alpha_2(\omega) \|x_1(t, \omega) - f(t, x_1(t, \omega))\|  \\
                        &\quad + \alpha_3(\omega) \|x_2(t, \omega) - f(t, x_2(t, \omega))\|  + \alpha_4(\omega) \|x_1(t, \omega) - f(t, x_2(t, \omega))\| \\
                        &\quad + \alpha_5(\omega) \|x_2(t, \omega) - f(t, x_1(t, \omega))\|
\endaligned
\eqno{(4.3)}
$$
for all $x_1(t, \omega), x_2(t, \omega)\in Q(\rho)$ and $\alpha_i : \Omega \rightarrow \mathbb{R}_+ \cup \{0\}$ for $i = 1, 2, 3, 4, 5$ such that $\sum^5_{i=1} \alpha_i(\omega) < 1$ almost surely;

 {\rm (3)}\, $h(t; \omega) \in D$,
\vskip 1mm

\noindent
then there exists a unique random solution of the equation $(4.1)$ in $Q(\rho)$ provided
\begin{eqnarray*}
  \|h(t, \omega)\|_D + l(\omega)\|f(t,0)\|_B\Big(\frac{1 + \alpha_3(\omega) + \alpha_4(\omega)}{1 - \alpha_2(\omega) + \alpha_5(\omega)}\Big)
  \leq \rho\Big(1 - \frac{l(\omega)}{1 - \alpha_2(\omega) - \alpha_5(\omega)}\Big),
\end{eqnarray*}
where the norm of $T(\omega)$ denoted by $l(\omega)$.
\end{thm}
\vskip 2mm

\begin{proof}
Let a mapping $\mathcal{U}(\omega) : Q(\rho) \rightarrow D$ defined by
\begin{eqnarray*}
 (\mathcal{U}x)(t, \omega) = h(t, \omega) + \int_S k(t, s, \omega)f(s, x(s, \omega))d_{\mu_0}(s).
\end{eqnarray*}
Then  we have
\begin{eqnarray*}
  \|(\mathcal{U}x)(t, \omega)\|_D &\leq& \|h(t, \omega)\|_D + l(\omega)\|f(t, x(t, \omega))\|_B \\ \nonumber
                        &\leq& \|h(t, \omega)\|_D + l(\omega)\|f(t, 0)\|_B + l(\omega)\|f(t, x(t, \omega)) - f(t, 0)\|_B.
\end{eqnarray*}
Thus it follows from (4.3) that
\begin{eqnarray*}
&&  \|f(t, x(t, \omega))-f(t, 0)\|_B \nonumber\\
&\leq& \alpha_1(\omega) \|x(t, \omega)\|_D + \alpha_2(\omega) \|x(t, \omega) - f(t, x(t, \omega))\|_D \\ \nonumber
                        && + \alpha_3(\omega) \|f(t, 0)\|_D + \alpha_4(\omega) \|x(t, \omega) - f(t, 0)\|_D  + \alpha_5(\omega) \|f(t, x(t, \omega))\|_D \\
                        &\leq& \alpha_1(\omega) \|x(t, \omega)\|_D + \alpha_2(\omega) \|x(t, \omega)\|_D \\
                        && + \alpha_2(\omega) \|f(t, x(t, \omega))-f(t, 0)\|_B + \alpha_2(\omega)\|f(t, 0)\|_D \\
                        && + \alpha_3(\omega)\|f(t, 0)\|_D + \alpha_4(\omega) \|x(t, \omega)\|_D + \alpha_4(\omega) \|f(t, 0)\|_D \\
                        && +  \alpha_5(\omega) \|f(t, x(t, \omega))-f(t, 0)\|_B + \alpha_5(\omega) \|f(t, 0)\|_D
\end{eqnarray*}
and so
\begin{eqnarray*}
 && (1 - \alpha_2(\omega) - \alpha_5(\omega))\|f(t, x(t, \omega))-f(t, 0)\|_B\\
  &\leq& (\alpha_1(\omega) + \alpha_2(\omega) + \alpha_4(\omega))\rho  + (\alpha_2(\omega) + \alpha_3(\omega) + \alpha_4(\omega) + \alpha_5(\omega))\|f(t,0)\|_D.
\end{eqnarray*}
Hence we have
$$
\aligned
  \|f(t, x(t, \omega))-f(t, 0)\|_B &\leq \Big(\frac{\alpha_1(\omega) + \alpha_2(\omega) + \alpha_4(\omega)}{1 - \alpha_2(\omega) - \alpha_5(\omega)}\Big)\rho \\
   &\quad + \Big(\frac{\alpha_2(\omega) + \alpha_3(\omega) + \alpha_4(\omega) + \alpha_5(\omega)}{1 - \alpha_2(\omega) - \alpha_5(\omega)}\Big)\|f(t,0)\|_D.
\endaligned
\eqno{(4.4)}
$$
Therefore, by (4.4), we have
$$
\aligned
&\quad  \|(\mathcal{U}x)(t, \omega)\|_D \\
&\leq \|h(t, \omega)\|_D + l(\omega)\|f(t, 0)\|_B  + \Big(\frac{\alpha_1(\omega) + \alpha_2(\omega) + \alpha_4(\omega)}{1 - \alpha_2(\omega) - \alpha_5(\omega)}\Big)l(\omega)\rho  \\
                        &\quad + \Big(\frac{\alpha_2(\omega) + \alpha_3(\omega) + \alpha_4(\omega) + \alpha_5(\omega)}{1 - \alpha_2(\omega) - \alpha_5(\omega)}\Big)l(\omega)\|f(t,0)\|_B \\
                        &\leq \|h(t, \omega)\|_D + \Big(\frac{\alpha_1(\omega) + \alpha_2(\omega) + \alpha_4(\omega)}{1 - \alpha_2(\omega) - \alpha_5(\omega)}\Big)l(\omega)\rho \\
                        &\quad + \Big(1 + \frac{\alpha_2(\omega) + \alpha_3(\omega) + \alpha_4(\omega) + \alpha_5(\omega)}{1 - \alpha_2(\omega) - \alpha_5(\omega)}\Big)l(\omega)\|f(t,0)\|_B \\
                        &\leq \|h(t, \omega)\|_D + \Big(\frac{\alpha_1(\omega) + \alpha_2(\omega) + \alpha_4(\omega)}{1 - \alpha_2(\omega) - \alpha_5(\omega)}\Big)l(\omega)\rho \\
                        &\quad + \Big(\frac{1 + \alpha_3(\omega) + \alpha_4(\omega)}{1 - \alpha_2(\omega) - \alpha_5(\omega)}\Big)l(\omega)\|f(t,0)\|_B \\
                        &< \rho
\endaligned
\eqno{(4.5)}
$$
and so, by (4.5), $(\mathcal{U}x)(t, \omega) \in Q(\rho)$.
Thus, for any $x_1(t, \omega)$, $x_2(t, \omega) \in Q(\rho)$ and, by the condition (2),  we have
\begin{eqnarray*}
 && \|(\mathcal{U}x_1)(t, \omega) - (\mathcal{U}x_2)(t, \omega)\|_D\nonumber \\
 &=& \Big\|\int_S k(t, s, \omega)[f(s, x_1(s, \omega)) - f(s, x_2(s, \omega))]d\mu_0(s)\Big\|_D \nonumber\\
        &\leq& l(\omega)\|f(s, x_1(s, \omega)) - f(s, x_2(s, \omega))\|_B  \nonumber\\
        &\leq& \alpha_1(\omega)\|x_1(t, \omega) - x_2(t, \omega)\|_D  + \alpha_2(\omega)\|x_1(t, \omega) - (\mathcal{U}x_1)(t, \omega)\|_D  \nonumber\\
        && + \alpha_3(\omega)\|x_2(t, \omega) - (\mathcal{U}x_2)(t, \omega)\|_D + \alpha_4(\omega)\|x_1(t, \omega) - (\mathcal{U}x_2)(t, \omega)\|_D  \nonumber\\
        && + \alpha_5(\omega)\|x_2(t, \omega) - (\mathcal{U}x_1)(t, \omega)\|_D.
\end{eqnarray*}
Consequently, $\mathcal{U}(\omega)$ is a random contractive mapping on $Q(\rho)$. Hence, by Theorem \ref{Maintheorem1}, there exists a random fixed point of $\mathcal{U}(\omega)$, which is the random solution of the equation (4.1). This completes the proof.
\end{proof}
\vskip 5mm

\textbf{Open Problem:} Can Theorems \ref{MainMotivate} and \ref{Maintheorem1} be generalized to non-separable Banach spaces?

\section*{{\bf Acknowledgements}}
The authors are gratefully thankful for referee's valuable comments, which significantly improve materials in this paper.
The first author was supported by Rajamangala University of Technology Lanna (RMUTL) for Ph.D. program at King Mongkut's University of Technology Thonburi (KMUTT). Yeol Je Cho was supported by Basic Science Research Program through the National Research Foundation of Korea (NRF) funded by the Ministry of Science, ICT and future Planning (2014R1A2A2A01002100).
This work was carried out while the third author (YJ. Cho) was visiting Theoretical and Computational Science Center (TaCS), Science Laboratory Building, Faculty of Science, King Mongkut's University of Technology Thonburi (KMUTT), Bangkok, Thailand, during 15 January-2 march, 2016. He thanks Professor Poom Kumam and the University for their hospitality and support.

 Moreover,  Poom Kumam was supported by the Thailand Research Fund (TRF) and the King Mongkut's University of Technology Thonburi (KMUTT) under the TRF Research Scholar Award  (Grant No. RSA6080047).



\vskip 5mm


\footnotesize

\begin{thebibliography}{10}

\bibitem{Achari}
Achari, J: On a pair of random generalized nonlinear contractions. Internat. J. Math. Math. Sci. 6, 467--475 (1983).

\bibitem{Arens}
Arens, RF: A topology for spaces of transformations, Ann.  Math. (2)47, 480-495  (1946).


\bibitem{Banach}
Banach, S: Sur les operations dans les ensembles abstraits et leur application aux equations integrals. Fund. Math. 3, 133--181 (1922).

\bibitem{BegShahzad}
Beg, I, Shahzad, N: Random fixed points of random multivalued operator on Polish spaces, Nonlinear Anal., 20 (1993), pp. 835--847.

\bibitem{BharuchaReid}
Bharucha-Reid, AT: Fixed point theorems in probabilistic analysis. Bull. Amer. Math. Soc. 82, 641--657 (1976).



\bibitem{ChoLiHuang}
Cho, YJ, Li, J, Huang, NJ: Random Ishikawa iterative sequence with errors for approximating random fixed points. Taiwanese Journal of Mathematics, Vol. 12, No. 1, pp. 51-61, February 2008.


\bibitem{Ciric01}
$\acute{\text{C}}$iri$\acute{\text{c}}$, CLj: On a generalization of a Gregus fixed point theorem. Czechoslov. Math. J. 50, 449--458 (2000).







\bibitem{Gregus01}
Gregus, M: A fixed point theorem in Banach space. Boll. Union. Mat. Ital. A 5(7), 193--198 (1980).


\bibitem{Hans1}
Hans, O: Random operator equations. In: Proceedings of 4th Berkeley Sympos. Math. Statist. Prob., vol. II, part I, pp. 185--202, University of California Press, Berkeley (1961).

\bibitem{Hans2}
Hanse, O: Reduzierende zufallige transformationen. Czechoslov. Math. J. 7(82), 154--158 (1957).



\bibitem{Itoh}
Itoh, S: Random fixed-point theorems with an application to random differential equations in Banach spaces. J. Math. Anal. Appl. 67, 261--273 (1979).



\bibitem{JoshiBose}
Joshi, MC, Bose, RK: Some Topics in Nonlinear Functional Analysis. Wiley, New York (1984).


\bibitem{JungChoKangLeeThakur}
Jung, JS,  Cho, YJ, Kang, SM, Lee, BS, Thakur, BS: Random fixed point theorems for a certain class of mappings in banach spaces. Czechoslovak Mathematical Journal, June 2000, Volume 50, Issue 2, pp 379-396.





\bibitem{Kumam02}
Kumam, P: Random common fixed points of single-valued and multivalued random operators in a uniformly convex
Banach space. J. Comput. Anal. Appl. 13, 368--375 (2011).

\bibitem{KumamKumam01}
Kumam, P, Kumam, W: Random fixed points of multivalued random operators with property (D). Random Oper.
Stoch. Equat. 15, 127--136 (2007).


\bibitem{KumamKumam02}
Kumam, W, Kumam, P: Random fixed point theorems for multivalued subsequentially limit-contractive maps
satisfying inwardness conditions. J. Comput. Anal. Appl. 14, 239--251 (2012).

\bibitem{KumamPlubtieng04}
Kumam, P, Plubtieng, S: Random fixed point theorems for asymptotically regular random operators. Demonst. Math. XLII, 131--141 (2009).



\bibitem{KumamPlubtieng08}
Kumam, P, Plubtieng, S: The characteristic of noncompact convexity and random fixed point theorem for set-valued operators. Czechoslov. Math. J. 57(132), 269--279 (2007).







\bibitem{KumamPlubtieng05}
Kumam, P, Plubtieng, S: Random fixed point theorems for multivalued nonexpansive non-self random operators. J.
Appl. Math. Stoch. Anal. 2006, Article ID 43796 (2006).






\bibitem{KumamPlubtieng01}
Kumam, P, Plubtieng, S: Random coincidence and random common fixed points of nonlinear multivalued random
operators. Thai J. Math. 5, 155--163 (2007) (Special issue).





\bibitem{LeePadgett}
Lee, ACH, Padgett, WJ: On random nonlinear contraction. Math. Syst. Theory 11, 77--84 (1977).





\bibitem{Mukherjee}
Mukherjee, A: Transformation aleatoires separable theorem all point fixed aleatoire. C.R. Acad. Sci. Paris, Ser. A-B 263, 393--395 (1966).



\bibitem{Padgett}
Padgett, WJ: On a nonlinear stochastic integral equation of the hammerstein type. Proc. Amer. Soc. 38, 625--631 (1973).


\bibitem{Radenovic}
Radenovi$\acute{\text{c}}$, S: Common fixed points under contractive conditions in cone metric spaces, Comput. Math. Appl. 58, 123--1278 (2009).

\bibitem{Rhoades}
Rhoades, BE: A comparison of various definitions of contractive mappings, Trans. Amer. Math. Soc. 226, 257-290 (1997).

\bibitem{Rothe}
Rothe, E: Zur theorie der topologischen ordnung und der Vektorfelder in Banachschen Raumen. Compos. Math. 5,
177--197 (1938).


\bibitem{Saha}
Saha, M: On some random fixed point of mappings over a Banach space with a probability measure. Proc. Natl. Acad. Sci. India, Sect. A 76, 219--224 (2006).

\bibitem{SahaDebnath}
Saha, M, Debnath, L: Random fixed point of mappings over a Hilbert space with a probability measure. Adv. Stud. Contemp. Math. 1, 79--84 (2007).


\bibitem{SahaDey}
Saha, M, Dey, D: Some random fixed point theorems for $(\theta, L)-$weak contractions. to appear in Hacet. J. Math. Stat.

\bibitem{SahaGanguly}
Saha, M, Ganguly, A: Random fixed point theorem on a $\acute{\text{C}}$iri$\acute{\text{c}}$-type contractive mapping and its consequence. Fixed Point Theory Appl. 2012, Article ID 209 (2012).


\bibitem{SehgalWaters}
Sehgal, VM, Waters, C: Some random fixed point theorems for condensing operators. Proc. Amer. Math. Soc. 90,
425--429 (1984).

\bibitem{Spacek}
Spacek, A: Zufallige Gleichungen. Czechoslov. Math. J. 5(80), 462--466 (1955).





\bibitem{Yosida}
Yosida, K: Functional analysis, Die Grundlehren der math. Wissenschaften,
Band 123, Academic Press, New York; Springer-Verlag, Berlin, 1965.




\end{thebibliography}
\end{document}